%% file: main.tex
\documentclass[11pt]{amsart}

\input{preamble}

\title[The multi-variable  $n/d$-conjecture]{%
  on the proof of multi-variable \\
  $n/d$-conjecture for hyperplane arrangements}

\author{Rui Tan}
\address{\resizebox{\dimexpr\linewidth-\parindent\relax}{!}{
  Rui Tan, School of Mathematical Sciences, Zhejiang University, Hangzhou, China}}
\email{ruitan@zju.edu.cn}

\author{Fanghan Xiang}
\address{\resizebox{\dimexpr\linewidth-\parindent\relax}{!}{
  Fanghan Xiang, School of Mathematical Sciences, Zhejiang University, Hangzhou, China}}
\email{fanghanxiang@zju.edu.cn}

\date{}

\begin{document}
\input{sections/00-abstract}
\maketitle

\input{sections/01-introduction}
\input{sections/02-bernstein-sato-nearby-cycles}
\input{sections/03-proof}

\input{references}
\end{document}

%% file: preamble.tex
\usepackage[T1]{fontenc}
\usepackage[utf8]{inputenc}
\usepackage{amsmath,amssymb,mathrsfs}
\usepackage{graphicx}
\usepackage{tikz-cd}
\usepackage{microtype}
\usepackage{needspace}
\usepackage{xcolor}
\definecolor{citationblue}{RGB}{0,70,180}
\usepackage[pagebackref]{hyperref}
\renewcommand*{\backref}[1]{}
\renewcommand*{\backrefalt}[4]{\ifnum#1>0 #2\fi}

\makeatletter
\g@addto@macro\@adjustvertspacing{%
  \setlength{\abovedisplayskip}{6pt}%
  \setlength{\belowdisplayskip}{6pt}%
  \setlength{\abovedisplayshortskip}{6pt}%
  \setlength{\belowdisplayshortskip}{6pt}%
}
\@adjustvertspacing

\let\ps@plain\ps@headings
\def\ps@firstpage{%
  \ps@empty
  \def\@oddhead{\normalfont\scriptsize\hfil\thepage}%
  \def\@evenhead{\normalfont\scriptsize\thepage\hfil}%
}
\makeatother
\theoremstyle{plain}
\newtheorem{theorem}{Theorem}[section]

\newtheorem{lemma}[theorem]{Lemma}

\newtheorem{conjecture}[theorem]{Conjecture}

\theoremstyle{definition}

\theoremstyle{remark}

\DeclareMathOperator{\Ann}{Ann}
\DeclareMathOperator{\Spec}{Spec}
\DeclareMathOperator{\Supp}{supp}

\hypersetup{
  colorlinks=true,
  linkcolor=blue,
  citecolor=green,
  urlcolor=citationblue
}

%% file: sections/00-abstract.tex
\begin{abstract}
In this article,
we use Wu's method to prove Budur's multi-variable $n/d$-conjecture for completeness.
\end{abstract}

%% file: sections/01-introduction.tex
\section{Introduction}

Let $X=\mathbb A^n_{\mathbb C}$ and let
$\mathbf f=(f_1,\ldots,f_r)$ be an $r$-tuple of nonzero polynomials on $X$.
Following the notation of \cite{Wu20},
we write
\[
  \mathbf s=(s_1,\ldots,s_r),\qquad
  \mathbf 1_r=(1,\ldots,1),\qquad
  \mathbf f^{\mathbf s}=\prod_{j=1}^r f_j^{s_j}.
\]
The Bernstein--Sato ideal of
$\mathbf f$ is
\[
  B_{\mathbf f}
  =\operatorname{Ann}_{\mathbb C[\mathbf s]}
  \left(
    \frac{\mathscr D_X[\mathbf s]\mathbf f^{\mathbf s}}
    {\mathscr D_X[\mathbf s]\mathbf f^{\mathbf s+\mathbf 1_r}}
  \right).
\]
We denote by $Z(B_{\mathbf f})$ the zero locus of $B_{\mathbf f}$.
For $r=1$, the monic generator of this ideal is the classical
Bernstein--Sato polynomial $b_f(s)$.
Budur proposed the multi-variable $n/d$-conjecture:

\begin{conjecture}[{\cite[Conjecture~3]{Bud13}}]\label{conj:budur}
Let $\mathbf f=(f_1,\ldots,f_r)$, where each $f_j$ defines a central
hyperplane arrangement in $\mathbb C^n$, not necessarily reduced, of
degree $d_j$, and $\prod_{j=1}^r f_j$ defines a central
essential indecomposable hyperplane arrangement. Then
\[
  \bigl\{\mathbf s\in\mathbb C^r:
    \textstyle\sum_{j=1}^r d_js_j+n=0\bigr\}
  \subseteq Z(B_{\mathbf f}).
\]
\end{conjecture}

The case $r=1$ in conjecture~\ref{conj:budur} is the $n/d$-conjecture of Budur, Musta\c{t}\u{a},
and Teitler\cite{BMT11}. 
Wu\cite{Wu26} and Davis-Yang\cite{DY26} both proved this conjecture,
more generally,
Davis-Yang\cite{DY26} proved Budur's multi-variable $n/d$-conjecture (conjecture~\ref{conj:budur}) by developing Sabbah's theory of multi-$V$-filtrations and using  specialization along the diagonal.

By using algebraic relative $\mathscr D$-modules and Beilinson--Bernstein
nearby cycles, Wu \cite{Wu26} proved that
\[
  b_f(-k/d)=0,\qquad k=n,n+1,\ldots,d,
\]
where $f$ defines a central essential indecomposable hyperplane
arrangement of degree $d$ in $\mathbb C^n$.
The case $k=n$ establishes the $n/d$-conjecture of Budur, Musta\c{t}\u{a},
and Teitler. 
Moreover,
it implies the strong monodromy conjecture for hyperplane arrangements\cite{BMT11}.

In this paper, 
we use Wu's method\cite{Wu26} to prove Budur's multi-variable $n/d$-conjecture for completeness.
More precisely,
we prove:

\begin{theorem}\label{thm:main}
Let $\mathbf f=(f_1,\ldots,f_r)$, where each $f_j$ defines a central
hyperplane arrangement in $\mathbb C^n$, not necessarily reduced, of
degree $d_j$, and $\prod_{j=1}^r f_j$ defines a central
essential indecomposable hyperplane arrangement. Set $d = \sum_{j=1}^{r}d_{j}$. Then, for every integer $k$ with $n\leq k\leq d$,
\[
  \bigl\{\mathbf s\in\mathbb C^r:
    \textstyle\sum_{j=1}^r d_js_j+k=0\bigr\}
  \subseteq Z(B_{\mathbf f}).
\]
In particular, Conjecture~\ref{conj:budur} holds.
\end{theorem}

%% file: sections/02-bernstein-sato-nearby-cycles.tex
\section{Multivariable Beilinson–Bernstein Nearby Cycles}\label{sec:nearby-cycles}
 
In this section,
we follow the notation in \cite[\S3]{Wu26}.
Let $X$ be a smooth algebraic variety over $\mathbb C$,
and let $\mathbf f=(f_1,\ldots,f_r)$ be an $r$-tuple of nonzero regular
functions on $X$. We write by $D$ the divisor $(\prod_{j=1}^{r}f_{j} = 0)$ and $j\colon U=X\setminus D\hookrightarrow X$ the open embedding. We denote by 
\[
  \mathbf f^{\mathbf s}=\prod_{j=1}^r f_j^{s_j}.
\]
and assume that $\prod_{j=1}^{r}f_{j}$ is not invertible. We write
\[
  \mathbb C[\mathbf s]=\mathbb C[s_1,\ldots,s_r],\qquad
  \mathscr D_X[\mathbf s]
  =\mathscr D_X\otimes_{\mathbb C}\mathbb C[\mathbf s].
\]
For a $\mathscr D_X[\mathbf s]$-module $\mathcal N$ and a prime ideal
$\mathfrak q\subseteq\mathbb C[\mathbf s]$, we have
\[
  \mathcal N_{\mathfrak q}
  =\mathcal N\otimes_{\mathbb C[\mathbf s]}
     \mathbb C[\mathbf s]_{\mathfrak q}
\]
and 
\[
  \Supp_{\mathbb C[\mathbf s]}(\mathcal N)
  =\{\mathfrak q\in\Spec\mathbb C[\mathbf s]:
        \mathcal N_{\mathfrak q}\neq0\}.
\]
Consider the left $\mathscr D_X[\mathbf s]$-module
\[
  \mathscr M
  =j_*\bigl(\mathcal O_U[\mathbf s]\mathbf f^{\mathbf s}\bigr).
\]
and $\mathscr M^r= \mathscr M \otimes_{\mathcal O_{X}} \omega_{X}$, which is the corresponding right $\mathscr D_X[\mathbf s]$-module under the side-change operation.

If $v$ is a global section of
\[
  \omega_X(*D):=\varinjlim_{k\to\infty}\omega_X(kD)
\]
such that $v \cdot \mathcal O_X|_{U}=\omega_U$, then we consider the coherent $\mathscr D_X[\mathbf s]$-submodule
\[ 
  \mathcal N:=\mathbf f^{\mathbf s} v \cdot \mathscr D_X[\mathbf s] \subseteq\mathscr M^r.
\]
For $\ell\in\mathbb Z$,
we define coherent $\mathscr D_X[\mathbf s]$-module
\[
  \mathcal N(\ell D)
  :=\mathbf f^{\mathbf s-\ell\mathbf1_r}v\cdot\mathscr D_X[\mathbf s].
\]
It is clear that 
\[\mathcal N(\ell D) \hookrightarrow \mathcal N((\ell +k) D)\]
for $k \in \mathbb{Z}_{\geq 0}$.

The Bernstein--Sato ideal associated with $\mathbf f$ and $v$ is defined by
\[
  B_{\mathbf f,v}
  :=\Ann_{\mathbb C[\mathbf s]}\!\left(
       \frac{\mathcal N}{\mathcal N(-D)}\right) .
\]

When $X=\mathbb A^n_{\mathbb C}=\Spec\mathbb C[z_1,\ldots,z_n]$
and $v=dz_1\wedge\cdots\wedge dz_n$, 
since the side-change operation commutes with 
$\mathbb C[\mathbf s]$-action, 
$B_{\mathbf f,v}$ coincides with the Bernstein--Sato ideal $B_{\mathbf f}$ of
$\mathbf f$.

\begin{lemma}\label{lem:parameter-support}
With the notation above,
\[
  Z(B_{\mathbf f,v})
  =\Supp_{\mathbb C[\mathbf s]}\!\left(
       \frac{\mathcal N}{\mathcal N(-D)}\right).
\]
Where $Z(B_{\mathbf f,v}) \subseteq\Spec\mathbb C[\mathbf s]$ is the zero locus of $B_{\mathbf f,v}$.
\end{lemma}

\begin{proof}
Take a finite affine open cover $X=\bigcup_i V_i$ such that
$\omega_{V_i}=\mathcal O_{V_i}e_i$ and $v|_{V_i}=h_i \cdot e_i$
for some rational function $h_i$ regular on $V_{i} \backslash D$.
By \cite[R\'esultat~2]{Mai23},
\[
  \frac{\mathscr D_{V_i}[\mathbf s]h_i\mathbf f^{\mathbf s}}
       {\mathscr D_{V_i}[\mathbf s]h_i
        \mathbf f^{\mathbf s+\mathbf1_r}}
\]
is algebraic relative holonomic.
Then the rest of the proof is similar to that of \cite[Lem.~3.1]{Wu26}.
\end{proof}

For a prime ideal $\mathfrak q\subseteq\mathbb C[\mathbf s]$.
We define
\[
  j_!\bigl(\mathscr M^r_{\mathfrak q}|_U\bigr)
  :=
  \mathbb D_Xj_*
  \bigl(\mathbb D_U(\mathscr M^r_{\mathfrak q}|_U)\bigr),
\]
where the definition of relative duality functor can be seen, for instance, in \cite[\S2]{Wu26}.

Since duality commutes with the side-change operation,
we obtain the following right module results from
\cite[Thms.~3.12--3.13]{Wu20}.

\begin{theorem}\label{thm:localized-extensions}
With the notation above, we have:
\begin{enumerate}
  \item
  $j_!(\mathscr M^r_{\mathfrak q}|_U)$ is concentrated in
  degree zero, i.e., $j_!(\mathscr M^r_{\mathfrak q}|_U)$ is a sheaf.

  \item
  The natural morphism
  $j_!(\mathscr M^r_{\mathfrak q}|_U) \longrightarrow j_*(\mathscr M^r_{\mathfrak q}|_U)$
  is injective.

  \item
  For $m\gg0$,
  $j_*(\mathscr M^r_{\mathfrak q}|_U)=\mathscr M^r_{\mathfrak q}=\mathcal N(mD)_{\mathfrak q}.$

  \item
  For $m\gg0$,
  $j_!(\mathscr M^r_{\mathfrak q}|_U)=\mathcal N(-mD)_{\mathfrak q}.$
\end{enumerate}
\end{theorem}

We define the multivariable Beilinson–Bernstein nearby cycle by
\begin{equation}\label{eq:nearby-definition}
  \Psi_{\mathbf f,\mathfrak q}
  :=
  \frac{j_*(\mathscr M^r_{\mathfrak q}|_U)}
       {j_!(\mathscr M^r_{\mathfrak q}|_U)}.
\end{equation}
By Theorem~\ref{thm:localized-extensions},
$\frac{\mathcal N(mD)_{\mathfrak q}}{\mathcal N((m-1)D)_{\mathfrak q}}$ is a subquotient
of $\Psi_{\mathbf f,\mathfrak q}$ for every $m \in \mathbb{Z}$.

%% file: sections/03-proof.tex
\section{Proof of the main result}\label{sec:proof}

Now let $X=\mathbb A^n_{\mathbb C}=\Spec\mathbb C[z_1,\ldots,z_n]$
and $v=dz_1\wedge\cdots\wedge dz_n$. Assume that $\mathbf f$ satisfies
the hypotheses of Conjecture~\ref{conj:budur}. We write 
\[
  d_j=\deg f_j,\qquad \mathbf d=(d_1,\ldots,d_r),\qquad
  d=\sum_{j=1}^r d_j.
\]
Fix an integer $n\leq k\leq d$ and write 
\[
  \mathfrak q=(\mathbf d\cdot\mathbf s+k)\subseteq\mathbb C[\mathbf s].
\]

\begin{theorem}\label{thm:first-shift-extension}
We have
\[
  j_!\bigl(\mathscr M^r_{\mathfrak q}|_U\bigr)
  =\mathcal N(-D)_{\mathfrak q}=\mathbf f^{\mathbf s+\mathbf1_r}v\cdot\mathscr D_X[\mathbf s]_{\mathfrak q}.
\]
\end{theorem}

\begin{proof}
By Theorem~\ref{thm:localized-extensions}(4),
it suffices to prove 
\[
  \frac{\mathcal N(-mD)_{\mathfrak q}
  }{\mathcal N(-(m+1)D)_{\mathfrak q}} = 0
\]
for all $m \in \mathbb{Z}_{> 0}$.
By \cite[Thm.~3.4(2)]{Wu20}, we can take
\[
  b(\mathbf s)=\prod_\nu(\mathbf a_\nu\cdot\mathbf s+c_\nu)
  \in B_{\mathbf f, v},
\]
 where $\mathbf a_\nu\in\mathbb Z_{\geq0}^r\setminus\{0\}$ and $c_\nu\in\mathbb Q_{>0}$.
Then by substitution,
$b(\mathbf s+m\mathbf1_r)$ annihilates
$\frac{\mathcal N(-mD)}{\mathcal N(-(m+1)D)}$.
Assume that $b(\mathbf s+m\mathbf1_r)\in\mathfrak q$.
Since $\mathfrak q$ is prime, we have $\mathbf a_\nu\cdot(\mathbf s+ m \mathbf{1}_r)+c_\nu = t(\mathbf d\cdot\mathbf s+k)$ for some $t$.
Comparing coefficients
\[
  \mathbf a_\nu=t\mathbf d
\]
hence $t>0$, 
then $c_\nu = t(k-md) \leq 0$, a contradiction.
Therefore,
$b(\mathbf s+m\mathbf1_r) \notin \mathfrak q$ and after localization at $\mathfrak q$,
 we obtain
\[
  \frac{\mathcal N(-mD)_{\mathfrak q}}
  {\mathcal N(-(m+1)D)_{\mathfrak q}} = 0.
\]
\end{proof}

Consequently,
\[
   \frac{\mathcal N_{\mathfrak q}}
  {\mathcal N(-D)_{\mathfrak q}} \hookrightarrow \Psi_{\mathbf f,\mathfrak q}.
\]
The image is generated by $[\mathbf f^{\mathbf s}v]$.
By Lemma~\ref{lem:parameter-support}, we have
\begin{equation}\label{eq:target-support}
  \{\mathbf d\cdot\mathbf s+k=0\}\subseteq Z(B_{\mathbf f, v})
  \quad\Longleftrightarrow\quad
  [\mathbf f^{\mathbf s}v]\neq0
  \quad\text{in }\Psi_{\mathbf f,\mathfrak q}.
\end{equation}

Now we prove Theorem~\ref{thm:main}. 
We follow the idea of \cite[\S4]{Wu26}.
Denote by $A$ the central essential indecomposable hyperplane arrangement defined by $\prod_{j=1}^{r}f_{j}$.
Let $L(A)$ be the intersection lattice of $A$, and let $\widetilde{L}(A)$ denote the set of all dense edges in $L(A)$.
For each $j=1,\ldots,r$, we may assume 
\[
  f_j=\prod_i f_{j,i}^{d_{j,i}},
\]
where each $f_{j,i}$ is a homogeneous linear polynomial and
$d_{j,i}\in\mathbb Z_{>0}$.
For every $W\in\widetilde{L}(A)$, define
\[
  d_{j,W}
  =\sum_{f_{j,i}|_W\equiv0}d_{j,i}
  \qquad\text{and}\qquad
  f_{j,W}
  =\prod_{f_{j,i}|_W\equiv0}f_{j,i}^{d_{j,i}}.
\]
We write
\[
  \mathbf d_W=(d_{1,W},\ldots,d_{r,W})
  \qquad\text{and}\qquad
  \mathbf f_W=(f_{1,W},\ldots,f_{r,W}),
\]
and set
\[
  d_W=\sum_{j=1}^r d_{j,W}
  \qquad\text{and}\qquad
  r_W=\operatorname{codim}_X W.
\]
Define
\[
  \Lambda_{\mathfrak q}
  =
  \left\{
    W\in\widetilde{L}(A)
    \ \middle|\
    (\mathbf d_W\cdot\mathbf s+\ell)=\mathfrak q
    \quad\text{for some }\ell\in\{r_W,\ldots,d_W\}
  \right\}.
\]
Since $\mathbf d_{\{0\}}=\mathbf d$, $d_{\{0\}}=d$,
and $r_{\{0\}}=n$, we have $\{0\}\in\Lambda_{\mathfrak q}$.
Thus,
$\Lambda_{\mathfrak q}$ is nonempty,
and then we may take  $W_0\in\Lambda_{\mathfrak q}$ such that
\[
  r_{W_0}
  =\min\{r_W | W\in\Lambda_{\mathfrak q}\},
\]
and an integer $\ell_0$ satisfying
\[
  r_{W_0}\leq\ell_0\leq d_{W_0}
  \qquad\text{and}\qquad
  (\mathbf d_{W_0}\cdot\mathbf s+\ell_0)=\mathfrak q.
\]
Since the global Bernstein--Sato ideal is always contained in the
local one,
we may replace $X$, $\mathbf f$, and $k$ by
$X/W_0$, $\mathbf f_{W_0}$, and $\ell_0$, respectively,
keeping the same notation.
Therefore,
we may assume $\Lambda_{\mathfrak q}=\{\{0\}\}$.

Let $\mu\colon Y\to X$ be the log resolution obtained by
successively blowing up the strict transforms of the unions
of the dense edges of dimension $m$, for $m=0,1,\ldots,n-2$.
Then
\[
  K_{Y/X}
  =\sum_{W\in\widetilde{L}(A)}(r_W-1)E_W
  \qquad\text{and}\qquad
  \mu^*D
  =\sum_{W\in\widetilde{L}(A)}d_WE_W,
\]
where $E_W$ is the irreducible divisor corresponding to $W$,
and $K_{Y/X}$ is the relative canonical divisor of $\mu$.
Set $\widetilde f_j=\mu^*f_j$ for $j=1,\ldots,r$,
$\widetilde{\mathbf f}=(\widetilde f_1,\ldots,\widetilde f_r)$,
and $\widetilde v=\mu^*v$, which is viewed as a section of $\omega_Y$
via $\mu^*\omega_X\hookrightarrow\omega_Y$.
Let
$\widetilde j\colon U\hookrightarrow Y$ be the open embedding.
Then we have a commutative diagram
\[
  \begin{tikzcd}
    U \arrow[r, "\widetilde j"] \arrow[d, "\mathrm{id}"']
      & Y \arrow[d, "\mu"] \\
    U \arrow[r, "j"] & X.
  \end{tikzcd}
\]
Since $j$ and $\widetilde j$ are affine open embeddings,
we have $j_+=j_*$ and $\widetilde j_+=\widetilde j_*$(see, for instance, \cite[\S2]{Wu26} for definition of $j_+$).
By functoriality, the above diagram yields
\begin{equation}\label{eq:resolution-star-extension}
  \mu_+\!\left(
    \widetilde j_*\bigl(\mathscr M^r_{\mathfrak q}|_U\bigr)
  \right)
  \simeq j_*\bigl(\mathscr M^r_{\mathfrak q}|_U\bigr).
\end{equation}
By \cite[Prop.~2.1]{Wu26}, we have
\begin{equation}\label{eq:resolution-shriek-extension}
  \mu_+\!\left(
    \widetilde j_!\bigl(\mathscr M^r_{\mathfrak q}|_U\bigr)
  \right)
  \simeq j_!\bigl(\mathscr M^r_{\mathfrak q}|_U\bigr).
\end{equation}
By Theorem~\ref{thm:localized-extensions}(2),
and the definition of multivariable nearby cycles, we have
a short exact sequence
\[
  0\longrightarrow
  \widetilde j_!\bigl(\mathscr M^r_{\mathfrak q}|_U\bigr)
  \longrightarrow
  \widetilde j_*\bigl(\mathscr M^r_{\mathfrak q}|_U\bigr)
  \longrightarrow
  \Psi_{\widetilde{\mathbf f},\mathfrak q}
  \longrightarrow0.
\]
Applying $\mu_+$ to the associated
distinguished triangle yields
\begin{equation}\label{eq:resolution-nearby-cycles}
  \mu_+\bigl(\Psi_{\widetilde{\mathbf f},\mathfrak q}\bigr)
  \simeq\Psi_{\mathbf f,\mathfrak q}.
\end{equation}
Now we consider
\[
  \widetilde{\mathcal N}
  :=\widetilde{\mathbf f}^{\mathbf s}\widetilde v
     \cdot\mathscr D_Y[\mathbf s]
  \subseteq\widetilde j_*\bigl(\mathscr M^r|_U\bigr).
\]
In a small neighborhood $V\subseteq Y$ of a point $p\in\bigcup_{W\in\widetilde L(A)}E_W$,
we have
\[
  \widetilde f_j|_V
  =u_j\prod_{W\in\widetilde L_p(A)}y_W^{d_{j,W}}
  \qquad\text{and}\qquad
  \widetilde v|_V
  =u'\prod_{W\in\widetilde L_p(A)}y_W^{r_W-1}\,dy,
\]
where $u_1,\ldots,u_r,u'$ are invertible regular functions on $V$,
$y_W$ is the local coordinate defining $E_W$,
$\widetilde L_p(A)
  :=\{W\in\widetilde L(A):p\in E_W\}$,
and $dy$ is the local volume form.
Moreover, we have 
\begin{equation}\label{eq:resolution-local-lattice}
  \widetilde{\mathcal N}|_V
  =u'\prod_{j=1}^r u_j^{s_j}
   \prod_{W\in\widetilde L_p(A)}
     y_W^{\mathbf d_W\cdot\mathbf s+r_W-1}
   \,dy\cdot\mathscr D_V[\mathbf s].
\end{equation}
Set  $g=u'\prod_{j=1}^r u_j^{s_j}$.
Since $u_1,\ldots,u_r,u'$ are invertible,
the conjugation action by $g$ induces an isomorphism of sheaves of rings
\[
   \mathscr D_V[\mathbf s]\simeq g^{-1} \cdot \mathscr D_V[\mathbf s] \cdot g,
\]
which sends each $\partial_{y_{W}}$ to 
\[
  g^{-1}(\partial_{y_W})g
  =\partial_{y_W}+\frac{\partial_{y_{W}} u'}{u'}
   +\sum_{j=1}^r s_j\frac{\partial_{y_W} u_j}{u_j}.
\]
Then we have an isomorphism 
\[
  \widetilde{\mathcal N}|_V \simeq \prod_{W\in\widetilde L_p(A)}
     y_W^{\mathbf d_W\cdot\mathbf s+r_W-1}
   \,dy\cdot\mathscr D_V[\mathbf s]
\]
after passing the $\mathscr D_V[\mathbf s]$-action by this conjugation,
and we may assume $g = 1$ for simplicity.

By Theorem~\ref{thm:localized-extensions}(4) and the local expression \eqref{eq:resolution-local-lattice}, locally around $p \in V$, we have
\[
  \widetilde j_!(\mathscr M^r_{\mathfrak q}|_U)|_V = \prod_{W\in\widetilde L_p(A)}
     y_W^{\mathbf d_W\cdot\mathbf s+r_W-1+md_W}
   \,dy\cdot\mathscr D_V[\mathbf s]_\mathfrak q
\]
for some $m\gg0$.
For $W\neq\{0\}$, we have
$\mathbf d_W\cdot\mathbf s+r_W -1 +\ell\notin\mathfrak q$
for every integer $\ell\geq 1$.
Otherwise, comparing coefficients
\[
  \mathbf{d}_W = t \mathbf{d} \qquad\text{and}\qquad  r_W-1+\ell=tk,
\]
hence
\[
  t = \frac{d_W}{d} \qquad\text{and}\qquad r_W \leq r_W-1+\ell=tk \leq d_W,
\]
contradicting $\Lambda_{\mathfrak q}=\{\{0\}\}$.
Similarly, for $W=\{0\}$,
$\mathbf d\cdot\mathbf s+ n-1 +\ell \notin \mathfrak q$ for every integer $\ell\neq k-n+1$.
Then 
\[
  \bigl(y_W^{\mathbf d_W\cdot\mathbf s+r_W-1+\ell}\,dy\bigr)
  \cdot\partial_{y_W}
  =-(\mathbf d_W\cdot\mathbf s+r_W-1+\ell)
    y_W^{\mathbf d_W\cdot\mathbf s+r_W-1+\ell-1}\,dy,
\]
yields 
\[
  \widetilde j_!\bigl(\mathscr M^r_{\mathfrak q}|_U\bigr)
  =
  \widetilde{\mathbf f}^{\mathbf s}
  \mathcal O_Y\bigl(-(k-n+1)E_{\{0\}}\bigr)
  \widetilde v\cdot\mathscr D_Y[\mathbf s]_{\mathfrak q}.
\]
We have a short exact sequence
\begin{equation}\label{eq:resolution-quotient-sequence}
  0\longrightarrow
  \mathcal Q_0:=
  \frac{\widetilde{\mathcal N}_{\mathfrak q}}
       {\widetilde j_!(\mathscr M^r_{\mathfrak q}|_U)}
  \longrightarrow
  \Psi_{\widetilde{\mathbf f},\mathfrak q}
  \longrightarrow
  \frac{\widetilde j_*(\mathscr M^r_{\mathfrak q}|_U)}
       {\widetilde{\mathcal N}_{\mathfrak q}}
  \longrightarrow0.
\end{equation}
By Theorem~\ref{thm:localized-extensions}(3) and
\eqref{eq:resolution-local-lattice}, locally around $p \in V$, 
$\frac{\widetilde j_*(\mathscr M^r_{\mathfrak q}|_U)}{\widetilde{\mathcal N}_{\mathfrak q}}$ has a finite filtration such that each successive quotient is of the form
\[
  \frac{
    y_{W'}^{-1}
    \displaystyle\prod_{W\in\widetilde L_p(A)}
      y_W^{\mathbf d_W\cdot\mathbf s+r_W-1-c_W}
    \,dy\cdot\mathscr D_V[\mathbf s]_{\mathfrak q}
  }{
    \displaystyle\prod_{W\in\widetilde L_p(A)}
      y_W^{\mathbf d_W\cdot\mathbf s+r_W-1-c_W}
    \,dy\cdot\mathscr D_V[\mathbf s]_{\mathfrak q},
  }
\]
where $c_W\in \mathbb{Z}_{\geq0}$ and $W'\in\widetilde L(A)$.
Such a quotient is nonzero if and only if
\[
  \mathbf d_{W'}\cdot\mathbf s+r_{W'}-1-c_{W'}
  \in\mathfrak q.
\]
After discarding the zero quotients, each remaining successive quotient $\mathcal Q_{W'}$ has the form
\begin{equation}\label{eq:filtration-quotient}
   \mathcal Q_{W'}|_V  \simeq \frac{
    y_{W'}^{-1}
    \displaystyle\prod_{W\in\widetilde L_p(A)}
      y_W^{\mathbf d_W\cdot\mathbf s+r_W-1-c_W}
    \,dy\cdot\mathscr D_V[\mathbf s]_{\mathfrak q}
  }{
    \displaystyle\prod_{W\in\widetilde L_p(A)}
      y_W^{\mathbf d_W\cdot\mathbf s+r_W-1-c_W}
    \,dy\cdot\mathscr D_V[\mathbf s]_{\mathfrak q}}
\end{equation}
for $W' \in  \Gamma_{\mathfrak q}$, where
\[
   \Gamma_{\mathfrak q}:=
  \left\{
    W\in\widetilde L(A)\ \middle|\
    \begin{array}{l}
      (\mathbf d_W\cdot\mathbf s+\ell)=\mathfrak q \quad 
      \text{for some }\ell\in \{1, \cdots, r_W-1\}
    \end{array}
  \right\}.
\]
In particular, $\{0\}\notin\Gamma_{\mathfrak q}$.
For $W'\in\Gamma_{\mathfrak q}$, necessarily
  $c_{W'}=r_{W'}-1-\frac{d_{W'}}{d}k$ and for other $c_W$ we can make suitable choices.
If $k=d$, then $\ell=d_W\geq r_W$,
so $\Gamma_{\mathfrak q}=\varnothing$ and
$\widetilde j_*(\mathscr M^r_{\mathfrak q}|_U)/\widetilde{\mathcal N}_{\mathfrak q}=0$.
Similarly, for $p\in E_{\{0\}}$, we have
\begin{equation}\label{eq:resolution-q0-local}
  \mathcal Q_0|_V\simeq
  \frac{
    y_{\{0\}}^{k-n}
    \displaystyle\prod_{\substack{W'\in\widetilde L(A)\\W'\ne\{0\}}}
      y_{W'}^{a_{W'}}
    \displaystyle\prod_{W\in\widetilde L_p(A)}
      y_W^{\mathbf d_W\cdot\mathbf s+r_W-1}
    \,dy\cdot\mathscr D_V[\mathbf s]_{\mathfrak q}
  }{
    y_{\{0\}}^{k-n+1}
    \displaystyle\prod_{\substack{W'\in\widetilde L(A)\\W'\ne\{0\}}}
      y_{W'}^{a_{W'}}
    \displaystyle\prod_{W\in\widetilde L_p(A)}
      y_W^{\mathbf d_W\cdot\mathbf s+r_W-1}
    \,dy\cdot\mathscr D_V[\mathbf s]_{\mathfrak q}
  },
\end{equation}
where we can choose any integers $a_{W'}\geq0$ since $\Lambda_{\mathfrak q}=\{\{0\}\}$.
Different choices of $a_{W'}$ mean we pick different local generators for $\mathcal Q_0$.

Locally around $p$, we consider $\mathcal Q_{W'}$ for every $W' \in \Gamma_{\mathfrak q}$.
By equation (\ref{eq:filtration-quotient}), 
$\mathfrak q \mathcal Q_{W'} = 0$ hence $\mathcal Q_{W'}$ is a $\mathscr D_Y \otimes_{\mathbb{C}} \kappa$-module,
where 
\[
  \kappa:=
  \mathbb C[\mathbf s]_{\mathfrak q}/
  \mathfrak q\mathbb C[\mathbf s]_{\mathfrak q}.
\]
We denote by $\bar{\mathbf s}$ the image of $\mathbf{s}$ in $\kappa$.
By relative Kashiwara's equivalence
(see, for instance, \cite[\S2.2]{GWX26}),
$\mathcal Q_{W'}$ corresponds to a $\mathscr{D}_{E_{W'}} \otimes_{\mathbb{C}} \kappa$-module generated by $\eta$ satisfying 
\[
  \eta \cdot(-\partial_{y_W}y_W)
  =(\mathbf d_W\cdot\bar{\mathbf s}+r_W-1-c_W)\eta
\]
for $W \neq W'$.
If $\{0\} \in \widetilde{L}_p(A)$,
then we have
\[
  \eta \cdot(-\partial_{y_{\{0\}}}y_{\{0\}})
  =(\mathbf{d} \cdot \bar{\mathbf{s}}+n-1-c_{\{0\}})\eta = (n-1-k-c_{\{0\}}) \eta.
\]
Moreover, for $j \geq 0$,
\[
  \eta y_{\{0\}}^{-j} \cdot (-\partial_{y_{\{0\}}}) = (n-1-k-c_{\{0\}}-j )\eta y_{\{0\}}^{-j-1}.
\]
Since $n-1-k-c_{\{0\}} < 0$,
$\mathcal Q_{W'}$ is an $\mathcal O_Y(*E_{\{0\}})$-module.
Then similar to the proof of \cite[Lem.~4.3]{Wu26},
we have
\begin{equation}\label{eq:resolution-quotient-vanishing}
  \mathcal H^i\mu_+\left(
    \frac{\widetilde j_*(\mathscr M^r_{\mathfrak q}|_U)}
         {\widetilde{\mathcal N}_{\mathfrak q}}
  \right)=0
  \qquad\text{for all }i<0.
\end{equation}

For simplicity, write
\[
  E=E_{\{0\}},
  \qquad
  E^\circ=E_{\{0\}}\setminus\bigcup_{W\ne\{0\}}E_W,
  \qquad
  B=E\setminus E^\circ.
\]
By equation (\ref{eq:resolution-q0-local}),
$\mathfrak q \mathcal Q_0=0$.
Thus $\mathcal Q_0$ is a $\mathscr D_Y\otimes_{\mathbb C}\kappa$-module,
by relative Kashiwara's equivalence, $\mathcal Q_0$ corresponds to a right
$\mathscr D_E\otimes_{\mathbb C}\kappa$-module $\mathcal E$ generated by the global section
\[ 
  e=
  \bigl[
    \mu^*(f_{1,1}^{k-n}) \cdot
    \widetilde{\mathbf f}^{\mathbf s}\widetilde v
  \bigr].
\] 
Locally, the section satisfies
\begin{equation}\label{eq:q0-logarithmic-action}
  e\cdot(-\partial_{y_W}y_W)
  =
  \bigl(\mathbf d_W\cdot\bar{\mathbf s} +r_W-1+a_W\bigr)e,
\end{equation}
where $W \neq \{0\}$ and  $a_W \geq 0$ is the multiplicity of $\mu^*(f_{1,1}^{k-n})$ along $y_W$.
We already have $\mathbf d_W\cdot\mathbf s+r_W -1 +\ell\notin\mathfrak q$
for every integer $\ell\geq 1$ and $W\neq\{0\}$.
Thus, 
its image  is nonzero in $\kappa$.
Then $a_W \geq 0$ yields 
\[
  \mathbf d_W\cdot\bar{\mathbf s}+r_W-1+a_W
  \notin\mathbb Z_{<0}.
\]

We consider the free module
\[
  \overline{\mathcal E}
  :=(\mathcal O_E\otimes_{\mathbb C}\kappa)e,
\]
it defines a right logarithmic connection on $(E,B)$ such that 
\[
    \mathcal E
    \simeq
    \overline{\mathcal E}
    \otimes_{\mathscr D_{E,B}\otimes_{\mathbb C}\kappa}
    (\mathscr D_E\otimes_{\mathbb C}\kappa)
    \simeq
    \overline{\mathcal E}
    \otimes_{\mathscr D_{E,B}}\mathscr D_E,
\]
where $\mathscr D_{E,B}$ denotes logarithmic differential operators
along $B$ (see~\cite{WZ21}).
By using $e$ to identify $\overline{\mathcal E}$ with $\mathcal O_E\otimes_{\mathbb C}\kappa$ and applying the standard logarithmic Spencer resolution construction with coefficient $\kappa$,
we obtain the logarithmic Spencer resolution of
$\overline{\mathcal E}$:
\[
     \cdots\longrightarrow{}
    K^{-r}
    \otimes_{\mathcal O_E}\mathscr D_{E,B}
    \longrightarrow{}
    K^{-r+1}
    \otimes_{\mathcal O_E}\mathscr D_{E,B}
    \longrightarrow{}\cdots 
    \longrightarrow{}
    \mathscr D_{E,B}\otimes_{\mathbb C}\kappa,
\]
where
\[ 
  K^{-r}
  =
  \bigwedge^r\mathcal T_E(-\log B)
  \otimes_{\mathbb C}\kappa,
  \qquad 0\le r\le n-1.
\]
Extending scalars in the above logarithmic Spencer resolution of $\overline{\mathcal E}$ to $\mathscr D_E$ yields a locally free right  $\mathscr D_E$-resolution of $\mathcal E$.
Therefore,
\[
  \mathcal E
  \otimes^{\mathbb L}_{\mathscr D_E}\mathcal O_E
\]
is represented by the complex $K^\bullet$,
whose differentials are determined by (\ref{eq:q0-logarithmic-action}) and the Leibniz rule.
We denote by $a \colon E \rightarrow \{0\}$,
then
\[
  a_{+}\mathcal E \simeq R\Gamma(E, K^{\bullet }),
\]
and the image of $e$ in $\mathcal{H}^{0}a_{+}\mathcal E$ is represented by $1 \in H^{0}(E,K^{0})$.
Consider the hypercohomology spectral sequence
\[
  E_1^{-r,j}=H^j(E,K^{-r})
  \Longrightarrow\mathbb H^{j-r}(E,K^\bullet).
\]
For each $\ell\ge1$, the only possible differential entering
$E_\ell^{0,0}$ is
\[
  d_\ell:E_\ell^{-\ell,\ell-1}
  \longrightarrow E_\ell^{0,0}.
\]
For $1 \leq \ell \leq n-1$,
by \cite[Lem.~4.4]{Wu26} 
\[
    E_{1}^{-\ell, \ell-1}=H^{\ell -1}(E,K^{-\ell})
    =H^{\ell-1}\left(
      E,\bigwedge^\ell\mathcal T_E(-\log B)
    \right)\otimes_{\mathbb C}\kappa=0.
\]
For $\ell \geq n$,
$E_{\ell}^{-\ell,\ell-1}$ lies outside the range of the complex.
Hence the source of $d_{\ell}$ vanishes.
Since $K^\bullet$ is concentrated in nonpositive degrees,
there are no differentials from $E_{\ell}^{0,0}$.
Therefore $E_1^{0,0}=E_\infty^{0,0}$,
and $E_{\infty}^{0,0}$ is a subspace of $\mathbb H^0(E,K^\bullet)$.
Hence we have 
\[
  H^0(E,K^0)= E_{1}^{0,0}
  \hookrightarrow\mathbb H^0(E,K^\bullet).
\]
Consequently,
the image of $e$ in $\mathcal{H}^{0}a_{+}\mathcal E$ is nonzero.
By Kashiwara's equivalence, its image in $\mathcal H^0\mu_+(\mathcal Q_0)$ is also nonzero.

Finally, by \eqref{eq:resolution-quotient-vanishing}, we have  an injection
\[
    \mathcal H^0\mu_+(\widetilde{\mathcal N}_{\mathfrak q})
    \hookrightarrow j_*(\mathscr M^r_{\mathfrak q}|_U),
\]
and consequently,
an injection
\[
  \mathcal H^0\mu_+(\mathcal Q_0)
    \hookrightarrow \Psi_{\mathbf f,\mathfrak q}.
\]
The global section $\widetilde{\mathbf f}^{\mathbf s}\widetilde v\otimes1$ of $\widetilde{\mathcal N}_{\mathfrak q} \otimes_{\mathcal O_{Y}} \mu^{*}\mathscr{D}_{X}$
induces morphisms in the derived category of $\mu^{-1}\mathcal O_X$-modules:
\[
  \mathcal O_Y \to \widetilde{\mathcal N}_{\mathfrak q}\otimes_{\mathcal O_Y}^{\mathbb L}\mu^*\mathscr D_X \to \widetilde{\mathcal N}_{\mathfrak q}\otimes_{\mathscr D_Y}^{\mathbb L}\mu^*\mathscr D_X \to \widetilde j_*(\mathscr M^r_{\mathfrak q}|_U)\otimes_{\mathscr D_Y}^{\mathbb L}\mu^*\mathscr D_X,
\]
where the first arrow sends $1$ to
$\widetilde{\mathbf f}^{\mathbf s}\widetilde v\otimes1$.
Applying $R\mu_*$ and composing with
$\mathcal O_X\to R\mu_*\mathcal O_Y$, we obtain
\[
  \mathcal O_X \to \mu_+(\widetilde{\mathcal N}_{\mathfrak q}) \to \mu_+\!\left(\widetilde j_*(\mathscr M^r_{\mathfrak q}|_U)\right) \simeq j_*(\mathscr M^r_{\mathfrak q}|_U).
\]
Then the global section $1$ of $\mathcal{O}_{X}$ yields a global section $w$ of $j_*(\mathscr M^r_{\mathfrak q}|_U)$,
which is the image of $\widetilde{\mathbf f}^{\mathbf s}\widetilde v$ under $\mu_{+}$.
Since $\mu|_U=\mathrm{id}_U$, we have
\[
  w|_U=\left.\widetilde{\mathbf f}^{\mathbf s}\widetilde v\right|_U=\left.\mathbf f^{\mathbf s}v\right|_U.
\]
Thus,
$w=\mathbf f^{\mathbf s}v$ since $j_*(\mathscr M^r_{\mathfrak q}|_U)$ has no nonzero sections supported on $D$,
and the image of $\widetilde{\mathbf f}^{\mathbf s}\widetilde v$ in $\mathcal H^0\mu_+(\widetilde{\mathcal N}_{\mathfrak q})$ is identified with $\mathbf f^{\mathbf s}v$ under the first injection.
Moreover,
the image of $\mu^{*}(f_{1,1}^{k-n})\cdot \widetilde{\mathbf f}^{\mathbf s}\widetilde v$ is identified with $f_{1,1}^{k-n} \cdot \mathbf f^{\mathbf s}v$.
We have shown that the image of 
\[ 
  e=
  \bigl[
    \mu^*(f_{1,1}^{k-n}) \cdot
    \widetilde{\mathbf f}^{\mathbf s}\widetilde v
  \bigr].
\] 
in $\mathcal H^0\mu_+(\mathcal Q_0)$ is nonzero.
Thus,
the image of $f_{1,1}^{k-n} \cdot \mathbf f^{\mathbf s} v$ is nonzero in $\Psi_{\mathbf f, \mathfrak q}$.
Hence $[\mathbf f^{\mathbf s}v]$ is nonzero in $\Psi_{\mathbf f, \mathfrak q}$, and \eqref{eq:target-support} completes the proof of Theorem~\ref{thm:main}.